\documentclass{article}%
\usepackage{amsmath}
\usepackage{amsfonts}
\usepackage{amssymb}
\usepackage{graphicx}%
\providecommand{\U}[1]{\protect\rule{.1in}{.1in}}
\date{1 March 2010}
\newtheorem{theorem}{Theorem}

\newtheorem{corollary}[theorem]{Corollary}

\newtheorem{definition}[theorem]{Definition}

\newtheorem{lemma}[theorem]{Lemma}

\newtheorem{remark}[theorem]{Remark}

\newenvironment{proof}[1][Proof]{\noindent\textbf{#1.} }{\ \rule{0.5em}{0.5em}}
\ifx\pdfoutput\relax\let\pdfoutput=\undefined\fi
\newcount\msipdfoutput
\ifx\pdfoutput\undefined\else
\ifcase\pdfoutput\else
\ifx\paperwidth\undefined\else
\ifdim\paperheight=0pt\relax\else\pdfpageheight\paperheight\fi
\ifdim\paperwidth=0pt\relax\else\pdfpagewidth\paperwidth\fi
\fi\fi\fi
\begin{document}

\title{Approximation by $q$-Sz\'{a}sz operators}
\author{N. I. Mahmudov\\Department of Mathematics\\Eastern Mediterranean University \\Gazimagusa, TRNC, Mersin 10 \\Turkey}
\maketitle

\begin{abstract}
This paper deals with approximating properties of the newly defined
$q$-generalization of the Sz\'{a}sz operators in the case $q>1$. Quantitative
estimates of the convergence in the polynomial weighted spaces and the
Voronovskaja's theorem are given. In particular, it is proved that the rate of
approximation by the $q$-Sz\'{a}sz operators ($q>1$) is of order $q^{-n}$
versus $1/n$ for the classical Sz\'{a}sz--Mirakjan operators.

\end{abstract}

\textbf{Keywords:} Positive linear operators, Sz\'{a}sz-Mirakjan operators,
Voronovskaja-type asymptotic formula, weighted space, direct results, Korovkin
theorem 

\textbf{Mathematics Subject Classification (2000):}  41A36, 41A25.

\section{Introduction}

The approximation of functions by using linear positive operators introduced
via $q$-Calculus is currently under intensive research. The pioneer work has
been made by A. Lupa\c{s} \cite{lupas} and G. M. Phillips \cite{phil} who
proposed generalizations of Bernstein polynomials based on the $q$-integers.
The $q$-Bernstein polynomials quickly gained the popularity, see
\cite{sof1}-\cite{vid}. Other important classes of discrete operators have
been investigated by using $q$-Calculus in the case $0<q<1$, for example
$q$-Meyer-K\"{o}nig operators \cite{trif}, \cite{dogru}, \cite{hep3},
$q$-Bleimann, Butzer and Hahn operators \cite{aral2}, \cite{mahsab},
\cite{mahsab2}, $q$-Sz\'{a}sz-Mirakjan operators \cite{gupta}, \cite{aral},
\cite{radu}, \cite{mah2}, $q$-Baskakov operators \cite{ag}.

In the present paper, we introduce a $q$-generalization of the Sz\'{a}sz
operators in the case $q>1$. Notice that different $q$-generalizations of
Sz\'{a}sz-Mirakjan operators were introduced and studied by A. Aral and V.
Gupta \cite{gupta}, \cite{aral}, by C. Radu \cite{radu} and by N. I. Mahmudov
\cite{mah2} in the case $0<q<1$. Since we define $q$-Sz\'{a}sz operators for
$q>1$, the rate of approximation by the $q$-Sz\'{a}sz operators ($q>1$) is of
order $q^{-n},$ which is essentially better than $1/n$ (rate of approximation
for the classical Sz\'{a}sz--Mirakjan operators). Thus our $q$-Sz\'{a}sz
operators have better approximation properties than the classical
Sz\'{a}sz--Mirakjan operators and the other $q$-Sz\'{a}sz-Mirakjan operators.

The paper is organized as follows. In Section 2, we give standard notations
that will be used throughout the paper, introduce $q$-Sz\'{a}sz operators and
evaluate the moments of $M_{n,q}$. In Section 3 we study convergence
properties of the $q$-Sz\'{a}sz operators in the polynomial weigthed spaces.
In Section 4, we give the quantitative Voronovskaja-type asymptotic formula.

\section{Construction of $M_{n,q}$ and estimation of moments}

Throughout the paper we employ the standard notations of $q$-calculus, see
\cite{kac}, \cite{ernst}.

$q$-integer and $q$-factorial are defined by%
\begin{align*}
\left[  n\right]  _{q}  &  :=\left\{
\begin{array}
[c]{c}%
\dfrac{1-q^{n}}{1-q}\ \ \ \ \text{if\ \ \ }q\in R^{+}\backslash\{1\},\\
n\ \ \ \ \ \ \ \text{if\ \ \ \ \ }q=1
\end{array}
\right.  \ \ \ \ \ \text{for }n\in N\ \ \ \ \text{and\ \ \ }\left[  0\right]
=0,\\
\left[  n\right]  _{q}!  &  :=\left[  1\right]  _{q}\left[  2\right]
_{q}...\left[  n\right]  _{q}\ \ \ \ \text{for }n\in N\ \ \ \ \text{and\ \ \ }%
\left[  0\right]  !=1.
\end{align*}
For integers $0\leq k\leq n$ $q$-binomial is defined by%
\[
\left[
\begin{array}
[c]{c}%
n\\
k
\end{array}
\right]  _{q}:=\frac{\left[  n\right]  _{q}!}{\left[  k\right]  _{q}!\left[
n-k\right]  _{q}!}.
\]
The $q$-derivative of a function $f(x)$, denoted by $D_{q}f$, is defined by%
\[
\left(  D_{q}f\right)  \left(  x\right)  :=\frac{f\left(  qx\right)  -f\left(
x\right)  }{\left(  q-1\right)  x},\ \ \ x\neq0,\ \ \ \left(  D_{q}f\right)
\left(  0\right)  :=\lim_{x\rightarrow0}\left(  D_{q}f\right)  \left(
x\right)  .
\]
The formula for the $q$-derivative of a product and quotient are%
\begin{equation}
D_{q}\left(  u\left(  x\right)  v\left(  x\right)  \right)  =D_{q}\left(
u\left(  x\right)  \right)  v\left(  x\right)  +u\left(  qx\right)
D_{q}\left(  v\left(  x\right)  \right)  . \label{q1}%
\end{equation}
Also, it is known that%
\begin{equation}
D_{q}x^{n}=\left[  n\right]  x^{n-1},\ \ \ D_{q}E\left(  ax\right)  =aE\left(
qax\right)  . \label{q2}%
\end{equation}
If $\left\vert q\right\vert >1$, or $\ 0<\left\vert q\right\vert <1$ and
$\left\vert z\right\vert <\frac{1}{1-q}$, the $q$-exponential function
$e_{q}\left(  x\right)  $ was defined by Jackson
\begin{equation}
e_{q}\left(  z\right)  :=\sum_{k=0}^{\infty}\frac{z^{k}}{\left[  k\right]
_{q}!}. \label{eq}%
\end{equation}
If $\left\vert q\right\vert >1$, $e_{q}\left(  z\right)  $ is an entire
function and%
\begin{equation}
e_{q}\left(  z\right)  =\prod_{j=0}^{\infty}\left(  1+\left(  q-1\right)
\frac{z}{q^{j+1}}\right)  ,\ \ \ \left\vert q\right\vert >1. \label{eq1}%
\end{equation}
There is another $q$-exponential function which is entire when $0<\left\vert
q\right\vert <1$ and which converges when $\left\vert z\right\vert <\frac
{1}{\left\vert 1-q\right\vert }$ if $\left\vert q\right\vert >1$. To obtain it
we must invert the base in (\ref{eq}), i.e. $q\rightarrow\frac{1}{q}$:%
\[
E_{q}\left(  z\right)  :=e_{1/q}\left(  z\right)  =\sum_{k=0}^{\infty}%
\frac{q^{k\left(  k-1\right)  /2}z^{k}}{\left[  k\right]  _{q}!}.
\]
We immediately obtain from (\ref{eq1}) that%
\[
E_{q}\left(  z\right)  =\prod_{j=0}^{\infty}\left(  1+\left(  1-q\right)
zq^{j}\right)  ,\ \ \ 0<\left\vert q\right\vert <1.
\]
The $q$-difference equations corresponding to $e_{q}\left(  z\right)  $ and
$E_{q}\left(  z\right)  $ are%
\begin{align*}
D_{q}e_{q}\left(  az\right)   &  =ae_{q}\left(  qz\right)  ,\ \ D_{q}%
E_{q}\left(  az\right)  =aE_{q}\left(  qaz\right)  ,\\
D_{1/q}e_{q}\left(  z\right)   &  =D_{1/q}E_{1/q}\left(  z\right)
=E_{1/q}\left(  q^{-1}z\right)  =e_{q}\left(  q^{-1}z\right)  ,\ \ \ \ q\neq0.
\end{align*}
Let $C_{p}$ is the set of all real valued functions $f$, continuous on
$\left[  0,\infty\right)  $ and such that $w_{p}f$ is uniformly continuous and
bounded on $\left[  0,\infty\right)  $ endowed with the norm%
\[
\left\Vert f\right\Vert _{p}:=\sup_{x\in\left[  0,\infty\right)  }w_{p}\left(
x\right)  \left\vert f\left(  x\right)  \right\vert .
\]
Here%
\[
w_{0}\left(  x\right)  :=1,\ \ \ \ w_{p}\left(  x\right)  :=\left(
1+x^{p}\right)  ^{-1},\ \ \ \ \text{if}\ \ \ p\in N.
\]
The corresponding Lipschitz classes are given for $0<\alpha\leq2$ by%
\begin{align*}
\Delta_{h}^{2}f\left(  x\right)   &  :=f\left(  x+2h\right)  -2f\left(
x+h\right)  +f\left(  x\right)  ,\\
\omega_{p}^{2}\left(  f;\delta\right)   &  :=\sup_{0<h\leq\delta}\left\Vert
\Delta_{h}^{2}f\right\Vert _{p},\ \ \ \ \ Lip_{p}^{2}\alpha:=\left\{  f\in
C_{p}:\omega_{p}^{2}\left(  f;\delta\right)  =0\left(  \delta^{\alpha}\right)
,\ \ \ \ \delta\rightarrow0^{+}\right\}  .
\end{align*}

Now we introduce the $q$-parametric Sz\'{a}sz operator.

\begin{definition}
Let $q>1$ and $n\in N$. For $f:\left[  0,\infty\right)  \rightarrow R$ we
define the Sz\'{a}sz operator based on the $q$-integers
\begin{equation}
M_{n,q}\left(  f;x\right)  :=\sum_{k=0}^{\infty}f\left(  \frac{\left[
k\right]  }{\left[  n\right]  }\right)  \frac{1}{q^{k\left(  k-1\right)  /2}%
}\frac{\left[  n\right]  ^{k}x^{k}}{\left[  k\right]  !}e_{q}\left(  -\left[
n\right]  q^{-k}x\right)  . \label{snq}%
\end{equation}

\end{definition}

Similarly as a classical Sz\'{a}sz operator $S_{n}$, the operator $M_{n,q}$ is
linear and positive. Furthermore, in the case of $q\rightarrow1^{+}$ we obtain
classical Sz\'{a}sz--Mirakjan operators.

Moments $M_{n,q}\left(  t^{m};x\right)  $ are of particular importance in the
theory of approximation by positive operators. From (\ref{snq}) we easily
derive the following recurrence formula and explicit formulas for moments
$M_{n,q}\left(  t^{m};x\right)  ,$ $m=0,1,2,3,4$.

\begin{lemma}
\label{rec}Let $q>1$. The following recurrence formula holds%
\begin{equation}
M_{n,q}\left(  t^{m+1};x\right)  =\sum_{j=0}^{m}\left(
\begin{array}
[c]{c}%
m\\
j
\end{array}
\right)  \frac{xq^{j}}{\left[  n\right]  ^{m-j}}M_{n,q}\left(  t^{j}%
;q^{-1}x\right)  . \label{rec1}%
\end{equation}

\end{lemma}

\begin{proof}
The recurrence formula (\ref{rec1}) easily follows from the definition of
$M_{n,q}$ and $q\left[  k\right]  +1=\left[  k+1\right]  .$
\begin{align*}
&  M_{n,q}\left(  t^{m+1};x\right) \\
&  =\sum_{k=0}^{\infty}\frac{\left[  k\right]  ^{m+1}}{\left[  n\right]
^{m+1}}\frac{1}{q^{k\left(  k-1\right)  /2}}\frac{\left[  n\right]  ^{k}x^{k}%
}{\left[  k\right]  !}e_{q}\left(  -\left[  n\right]  q^{-k}x\right) \\
&  =\sum_{k=1}^{\infty}\frac{\left[  k\right]  ^{m}}{\left[  n\right]  ^{m}%
}\frac{1}{q^{k\left(  k-1\right)  /2}}\frac{\left[  n\right]  ^{k-1}x^{k}%
}{\left[  k-1\right]  !}e_{q}\left(  -\left[  n\right]  q^{-k}x\right) \\
&  =\sum_{k=0}^{\infty}\frac{\left(  q\left[  k\right]  +1\right)  ^{m}%
}{\left[  n\right]  ^{m}}\frac{1}{q^{k\left(  k+1\right)  /2}}\frac{\left[
n\right]  ^{k}x^{k+1}}{\left[  k\right]  !}e_{q}\left(  -\left[  n\right]
q^{-k}q^{-1}x\right) \\
&  =\sum_{k=0}^{\infty}\frac{1}{\left[  n\right]  ^{m}}\sum_{j=0}^{m}\left(
\begin{array}
[c]{c}%
m\\
j
\end{array}
\right)  q^{j}\left[  k\right]  ^{j}\frac{1}{q^{k\left(  k+1\right)  /2}}%
\frac{\left[  n\right]  ^{k}x^{k+1}}{\left[  k\right]  !}e_{q}\left(  -\left[
n\right]  q^{-k}q^{-1}x\right) \\
&  =\sum_{j=0}^{m}\left(
\begin{array}
[c]{c}%
m\\
j
\end{array}
\right)  \frac{xq^{j}}{\left[  n\right]  ^{m-j}}\sum_{k=0}^{\infty}%
\frac{\left[  k\right]  ^{j}}{\left[  n\right]  ^{j}}\frac{1}{q^{k\left(
k-1\right)  /2}}\frac{\left[  n\right]  ^{k}x^{k}}{\left[  k\right]  !q^{k}%
}e_{q}\left(  -\left[  n\right]  q^{-k}q^{-1}x\right) \\
&  =\sum_{j=0}^{m}\left(
\begin{array}
[c]{c}%
m\\
j
\end{array}
\right)  \frac{xq^{j}}{\left[  n\right]  ^{m-j}}M_{n,q}\left(  t^{j}%
;q^{-1}x\right)  .
\end{align*}

\end{proof}

\begin{lemma}
\label{lemq}The following identities hold for all $q>1,$ $x\in\left[
0,\infty\right)  $, $n\in N$, and $k\geq0:$%
\begin{align}
xD_{q}s_{nk}\left(  q;x\right)   &  =\left[  n\right]  \left(  \frac{\left[
k\right]  }{\left[  n\right]  }-x\right)  s_{nk}\left(  q;x\right)
,\nonumber\\
M_{n,q}\left(  t^{m+1};x\right)   &  =\frac{x}{\left[  n\right]  }D_{q}%
M_{n,q}\left(  t^{m};x\right)  +xM_{n,q}\left(  t^{m};x\right)  . \label{rec3}%
\end{align}

\end{lemma}

\begin{proof}
The \ first identitiy follows from the following simple calculations%
\begin{align*}
xD_{q}s_{nk}\left(  q;x\right)   &  =\left[  k\right]  _{q}\frac
{1}{q^{k\left(  k-1\right)  /2}}\frac{\left[  n\right]  _{q}^{k}x^{k}}{\left[
k\right]  _{q}!}e_{q}\left(  -\left[  n\right]  q^{-k}x\right)  -xq^{-k}%
\left[  n\right]  _{q}\frac{1}{q^{k\left(  k-1\right)  /2}}\frac{\left[
n\right]  _{q}^{k}q^{k}x^{k}}{\left[  k\right]  _{q}!}e_{q}\left(  -\left[
n\right]  _{q}q^{-k}x\right) \\
&  =\left[  k\right]  _{q}s_{nk}\left(  q;x\right)  -x\left[  n\right]
_{q}s_{nk}\left(  q;x\right)  =\left[  n\right]  \left(  \frac{\left[
k\right]  }{\left[  n\right]  }-x\right)  s_{nk}\left(  q;x\right)  .
\end{align*}
The second one follows from the first.%
\begin{align*}
xD_{q}M_{n,q}\left(  t^{m};x\right)   &  =\left[  n\right]  \sum_{k=0}%
^{\infty}\left(  \frac{\left[  k\right]  }{\left[  n\right]  }\right)
^{m}\left(  \frac{\left[  k\right]  }{\left[  n\right]  }-x\right)
s_{nk}\left(  q;x\right) \\
&  =\left[  n\right]  \sum_{k=0}^{\infty}\left(  \frac{\left[  k\right]
}{\left[  n\right]  }\right)  ^{m+1}s_{nk}\left(  q;x\right)  -\left[
n\right]  x\sum_{k=0}^{\infty}\left(  \frac{\left[  k\right]  }{\left[
n\right]  }\right)  ^{m}s_{nk}\left(  q;x\right) \\
&  =\left[  n\right]  M_{n,q}\left(  t^{m+1};x\right)  -\left[  n\right]
xM_{n,q}\left(  t^{m};x\right)  .
\end{align*}

\end{proof}

\begin{lemma}
\label{moments}Let $q>1$. We have%
\begin{align*}
M_{n,q}\left(  1;x\right)   &  =1,\ \ \ M_{n,q}\left(  t;x\right)
=x,\ \ \ M_{n,q}\left(  t^{2};x\right)  =x^{2}+\frac{1}{\left[  n\right]
}x,\\
M_{n,q}\left(  t^{3};x\right)   &  =x^{3}+\frac{2+q}{\left[  n\right]  }%
x^{2}+\frac{1}{\left[  n\right]  ^{2}}x,\\
M_{n,q}\left(  t^{4};x\right)   &  =x^{4}+\left(  3+2q+q^{2}\right)
\frac{x^{3}}{\left[  n\right]  }+\left(  3+3q+q^{2}\right)  \frac{x^{2}%
}{\left[  n\right]  ^{2}}+\frac{1}{\left[  n\right]  ^{3}}x.
\end{align*}

\end{lemma}

\begin{proof}
For a fixed $x\in R_{+}$, by the $q$-Taylor theorem \cite{kac}, we obtain%
\[
\varphi_{n}\left(  t\right)  =\sum_{k=0}^{\infty}\frac{\left(  t-x\right)
_{1/q}^{k}}{\left[  k\right]  _{1/q}!}D_{1/q}^{k}\varphi_{n}\left(  x\right)
.
\]
Choosing $t=0$ and taking into account%
\[
\left(  -x\right)  _{1/q}^{k}=\left(  -1\right)  ^{k}x^{k}q^{-k\left(
k-1\right)  /2},\ \ \ \ D_{1/q}^{k}e_{q}\left(  -\left[  n\right]
_{q}x\right)  =\left(  -1\right)  ^{k}q^{-k\left(  k-1\right)  /2}\left[
n\right]  _{q}^{k}e_{q}\left(  -\left[  n\right]  _{q}q^{-k}x\right)
\]
we get for $\varphi_{n}\left(  x\right)  =e_{q}\left(  -\left[  n\right]
x\right)  $ that
\begin{align*}
1  &  =\varphi_{n}\left(  0\right)  =\sum_{k=0}^{\infty}\frac{\left(
-1\right)  ^{k}x^{k}}{q^{k\left(  k-1\right)  /2}\left[  k\right]  _{1/q}%
!}D_{1/q}^{k}\varphi_{n}\left(  x\right) \\
&  =\sum_{k=0}^{\infty}\frac{\left(  -1\right)  ^{k}x^{k}}{\left[  k\right]
_{q}!}\left(  -1\right)  ^{k}q^{-k\left(  k-1\right)  /2}\left[  n\right]
_{q}^{k}e_{q}\left(  -\left[  n\right]  _{q}q^{-k}x\right) \\
&  =\sum_{k=0}^{\infty}\frac{\left[  n\right]  ^{k}x^{k}}{\left[  k\right]
!q^{k\left(  k-1\right)  /2}}e_{q}\left(  -\left[  n\right]  q^{-k}z\right)  .
\end{align*}
In other words $M_{n,q}\left(  1;x\right)  =1$.

Calculation of $M_{n,q}\left(  t^{i};x\right)  $, $i=1,2,3,4,$ based on the
recurrence formula (\ref{rec3}) ( or (\ref{rec1})). We only calculate
$M_{n,q}\left(  t^{3};x\right)  $ and $M_{n,q}\left(  t^{4};x\right)  :$%
\begin{align*}
M_{n,q}\left(  t^{3};x\right)   &  =\frac{x}{\left[  n\right]  }D_{q}%
M_{n,q}\left(  t^{2};x\right)  +xM_{n,q}\left(  t^{2};x\right) \\
&  =\frac{x}{\left[  n\right]  }\left(  \left[  2\right]  x+\frac{1}{\left[
n\right]  }\right)  +x\left(  x^{2}+\frac{1}{\left[  n\right]  }x\right) \\
&  =\frac{1}{\left[  n\right]  ^{2}}x+\frac{2+q}{\left[  n\right]  }%
x^{2}+x^{3}.
\end{align*}%
\begin{align*}
M_{n,q}\left(  t^{4};x\right)   &  =\frac{x}{\left[  n\right]  }D_{q}%
M_{n,q}\left(  t^{3};x\right)  +xM_{n,q}\left(  t^{3};x\right) \\
&  =\frac{x}{\left[  n\right]  }\left(  \frac{1}{\left[  n\right]  ^{2}}%
+\frac{2+q}{\left[  n\right]  }\left[  2\right]  x+\left[  3\right]
x^{2}\right)  +x\left(  \frac{1}{\left[  n\right]  ^{2}}x+\frac{2+q}{\left[
n\right]  }x^{2}+x^{3}\right) \\
&  =\frac{1}{\left[  n\right]  ^{3}}x+\left(  3+3q+q^{2}\right)  \frac{x^{2}%
}{\left[  n\right]  ^{2}}+\left(  3+2q+q^{2}\right)  \frac{x^{3}}{\left[
n\right]  }+x^{4}.
\end{align*}

\end{proof}

\begin{lemma}
\label{lema}Assume that $q>1$. For every $x\in\left[  0,\infty\right)  $ there
hold%
\begin{align}
M_{n,q}\left(  \left(  t-x\right)  ^{2};x\right)   &  =\frac{x}{\left[
n\right]  },\label{m2}\\
M_{n,q}\left(  \left(  t-x\right)  ^{3};x\right)   &  =\frac{1}{\left[
n\right]  ^{2}}x+\left(  q-1\right)  \frac{x^{2}}{\left[  n\right]  },\\
M_{n,q}\left(  \left(  t-x\right)  ^{4};x\right)   &  =\frac{1}{\left[
n\right]  ^{3}}x+\left(  q^{2}+3q-1\right)  \frac{x^{2}}{\left[  n\right]
^{2}}+\left(  q-1\right)  ^{2}\frac{x^{3}}{\left[  n\right]  }\text{.}
\label{m3}%
\end{align}

\end{lemma}

\begin{proof}
First of all we give an explicit formula for $M_{n,q}\left(  \left(
t-x\right)  ^{4};x\right)  .$%
\begin{align*}
M_{n,q}\left(  \left(  t-x\right)  ^{3};x\right)   &  =M_{n,q}\left(
t^{3};x\right)  -3xM_{n,q}\left(  t^{2};x\right)  +3x^{2}M_{n,q}\left(
t;x\right)  -x^{3}\\
&  =x^{3}+\frac{2+q}{\left[  n\right]  }x^{2}+\frac{1}{\left[  n\right]  ^{2}%
}x-3x\left(  x^{2}+\frac{x}{\left[  n\right]  }\right)  +3x^{3}-x^{3}\\
&  =\frac{1}{\left[  n\right]  ^{2}}x+\left(  q-1\right)  \frac{x^{2}}{\left[
n\right]  }.
\end{align*}%
\begin{align*}
M_{n,q}\left(  \left(  t-x\right)  ^{4};x\right)   &  =M_{n,q}\left(
t^{4};x\right)  -4xM_{n,q}\left(  t^{3};x\right)  +6x^{2}M_{n,q}\left(
t^{2};x\right)  -4x^{3}M_{n,q}\left(  t;x\right)  +x^{4}\\
&  =\frac{1}{\left[  n\right]  ^{3}}x+\left(  3+3q+q^{2}\right)  \frac{x^{2}%
}{\left[  n\right]  ^{2}}+\left(  3+2q+q^{2}\right)  \frac{x^{3}}{\left[
n\right]  }+x^{4}\\
&  -4x\left(  \frac{1}{\left[  n\right]  ^{2}}x+\frac{2+q}{\left[  n\right]
}x^{2}+x^{3}\right)  +6x^{2}\left(  x^{2}+\frac{x}{\left[  n\right]  }\right)
-4x^{4}+x^{4}\\
&  =\frac{1}{\left[  n\right]  ^{3}}x+\left(  -1+3q+q^{2}\right)  \frac{x^{2}%
}{\left[  n\right]  ^{2}}+\left(  q-1\right)  ^{2}\frac{x^{3}}{\left[
n\right]  }.
\end{align*}

\end{proof}

Now we prove explicit formula for the moments $M_{n,q}\left(  t^{m};x\right)
$, which a $q$-analogue of a result of Becker, see \cite{becker} Lemma 3.

\begin{lemma}
\label{lemfor}For $q>1,$ $m\in N$ there holds%
\begin{equation}
M_{n,q}\left(  t^{m};x\right)  =\sum_{j=1}^{m}\mathbb{S}_{q}\left(
m,j\right)  \frac{x^{j}}{\left[  n\right]  ^{m-j}}, \label{rec2}%
\end{equation}
where
\begin{align}
\mathbb{S}_{q}\left(  m+1,j\right)   &  =\left[  j\right]  \mathbb{S}%
_{q}\left(  m,j\right)  +\mathbb{S}_{q}\left(  m,j-1\right)  ,\ \ m\geq
0,\ j\geq1,\nonumber\\
\mathbb{S}_{q}\left(  0,0\right)   &  =1,\ \ \mathbb{S}_{q}\left(  m,0\right)
=0,\ \ m>0,\ \ \mathbb{S}_{q}\left(  m,j\right)  =0,\ \ \ m<j. \label{str}%
\end{align}
In particular $M_{n,q}\left(  t^{m};x\right)  $ is a polynomial of degree $m$
without a constant term.
\end{lemma}

\begin{proof}
Because of $M_{n,q}\left(  t;x\right)  =x$, $M_{n,q}\left(  t^{2};x\right)
=x^{2}+\dfrac{x}{\left[  n\right]  }$, the representation (\ref{rec2}) holds
true for $m=1,2$ with $\mathbb{S}_{q}\left(  2,1\right)  =1$, $\mathbb{S}%
_{q}\left(  1,1\right)  =1$.

Now assume (\ref{rec2}) to be valued for $m$ then by Lemma \ref{lemq} we have
\begin{align*}
M_{n,q}\left(  t^{m+1};x\right)   &  =\frac{x}{\left[  n\right]  }D_{q}%
M_{n,q}\left(  t^{m};x\right)  +xM_{n,q}\left(  t^{m};x\right) \\
&  =\frac{x}{\left[  n\right]  }\sum_{j=1}^{m}\left[  j\right]  \mathbb{S}%
_{q}\left(  m,j\right)  \frac{x^{j-1}}{\left[  n\right]  ^{m-j}}+x\sum
_{j=1}^{m}\mathbb{S}_{q}\left(  m,j\right)  \frac{x^{j}}{\left[  n\right]
^{m-j}}\\
&  =\sum_{j=1}^{m}\left[  j\right]  \mathbb{S}_{q}\left(  m,j\right)
\frac{x^{j}}{\left[  n\right]  ^{m-j+1}}+\sum_{j=1}^{m}\mathbb{S}_{q}\left(
m,j\right)  \frac{x^{j+1}}{\left[  n\right]  ^{m-j}}\\
&  =\frac{x}{\left[  n\right]  ^{m}}\mathbb{S}_{q}\left(  m,1\right)
+x^{m+1}\mathbb{S}_{q}\left(  m,m\right) \\
&  +\sum_{j=2}^{m}\left(  \left[  j\right]  \mathbb{S}_{q}\left(  m,j\right)
+\mathbb{S}_{q}\left(  m,j-1\right)  \right)  \frac{x^{j}}{\left[  n\right]
^{m-j+1}}.
\end{align*}

\end{proof}

\begin{remark}
For $q=1$ the formulae (\ref{str}) become recurrence formulas satisfied by
Stirling numbers of the second type.
\end{remark}

\section{ $M_{n,q}$ in polynomial weighted spaces}

\begin{lemma}
\label{lemb}Let $p\in N\cup\left\{  0\right\}  $ and $q\in\left(
1,\infty\right)  $ be fixed. Then there exists a positive constant
$K_{1}\left(  q,p\right)  $ such that%
\begin{equation}
\left\Vert M_{n,q}\left(  1/w_{p};x\right)  \right\Vert _{p}\leq K_{1}\left(
q,p\right)  ,\ \ \ \ n\in N. \label{b1}%
\end{equation}
Moreover for every $f\in C_{p}$ we have
\begin{equation}
\left\Vert M_{n,q}\left(  f\right)  \right\Vert _{p}\leq K_{1}\left(
q,p\right)  \left\Vert f\right\Vert _{p},\ \ \ \ n\in N. \label{b2}%
\end{equation}
Thus $M_{n,q}$ is a linear positive operator from $C_{p}$ into $C_{p}$ for any
$p\in N\cup\left\{  0\right\}  $.
\end{lemma}

\begin{proof}
The inequality (\ref{b1}) is obvious for $p=0$. Let $p\geq1.$ Then by
(\ref{rec2}) we have%
\[
w_{p}\left(  x\right)  M_{n,q}\left(  1/w_{p};x\right)  =w_{p}\left(
x\right)  +w_{p}\left(  x\right)  \sum_{j=1}^{p}\mathbb{S}_{q}\left(
p,j\right)  \frac{x^{j}}{\left[  n\right]  ^{p-j}}\leq K_{1}\left(
q,p\right)  ,
\]
$K_{1}\left(  q,p\right)  $ is a positive constant depending on $p$ and $q$.
From this follows (\ref{b1}). On the other hand
\[
\left\Vert M_{n,q}\left(  f\right)  \right\Vert _{p}\leq\left\Vert
f\right\Vert _{p}\left\Vert M_{n,q}\left(  1/w_{p}\right)  \right\Vert _{p}%
\]
for every $f\in C_{p}$. By applying (\ref{b1}), we obtain (\ref{b2}).
\end{proof}

\begin{lemma}
Let $p\in N\cup\left\{  0\right\}  $ and $q\in\left(  1,\infty\right)  $ be
fixed. Then there exists a positive constant $K_{2}\left(  q,p\right)  $ such
that%
\begin{equation}
\left\Vert M_{n,q}\left(  \frac{\left(  t-\cdot\right)  ^{2}}{w_{p}\left(
t\right)  };\cdot\right)  \right\Vert _{p}\leq\frac{K_{2}\left(  q,p\right)
}{\left[  n\right]  },\ \ \ \ n\in N. \label{k2}%
\end{equation}

\end{lemma}

\begin{proof}
The formula (\ref{m2}) imply (\ref{k2}) for $p=0$. We have%
\[
M_{n,q}\left(  \frac{\left(  t-x\right)  ^{2}}{w_{p}\left(  t\right)
};x\right)  =M_{n,q}\left(  \left(  t-x\right)  ^{2};x\right)  +M_{n,q}\left(
\left(  t-x\right)  ^{2}t^{p};x\right)  ,
\]
for $p,n\in N$. If $p=1$ then we get%
\begin{align*}
M_{n,q}\left(  \left(  t-x\right)  ^{2}\left(  1+t\right)  ;x\right)   &
=M_{n,q}\left(  \left(  t-x\right)  ^{2};x\right)  +M_{n,q}\left(  \left(
t-x\right)  ^{2}t;x\right) \\
&  =M_{n,q}\left(  \left(  t-x\right)  ^{3};x\right)  +\left(  1+x\right)
M_{n,q}\left(  \left(  t-x\right)  ^{2};x\right)  ,
\end{align*}
which by Lemma \ref{lema} yields (\ref{k2}) for $p=1$.

Let $p\geq2$. By applying (\ref{rec2}), we get%
\begin{align*}
&  w_{p}\left(  x\right)  M_{n,q}\left(  \left(  t-x\right)  ^{2}%
t^{p};x\right) \\
&  =w_{p}\left(  x\right)  \left(  M_{n,q}\left(  t^{p+2};x\right)
-2xM_{n,q}\left(  t^{p+1};x\right)  +x^{2}M_{n,q}\left(  t^{p};x\right)
\right) \\
&  =w_{p}\left(  x\right)  \left(  x^{p+2}+\sum_{j=1}^{p+1}\mathbb{S}%
_{q}\left(  p+2,j\right)  \frac{x^{j}}{\left[  n\right]  ^{p+2-j}}%
-2x^{p+2}-2\sum_{j=1}^{p}\mathbb{S}_{q}\left(  p+1,j\right)  \frac{x^{j+1}%
}{\left[  n\right]  ^{p+1-j}}+x^{p+2}+\sum_{j=1}^{p-1}\mathbb{S}_{q}\left(
p,j\right)  \frac{x^{j+2}}{\left[  n\right]  ^{p-j}}\right) \\
&  =w_{p}\left(  x\right)  \left(  \sum_{j=2}^{p}\left(  \mathbb{S}_{q}\left(
p+2,j\right)  -2\mathbb{S}_{q}\left(  p+1,j\right)  +\mathbb{S}_{q}\left(
p,j\right)  \right)  \frac{x^{j+1}}{\left[  n\right]  ^{p+1-j}}\right. \\
&  +\left.  \mathbb{S}_{q}\left(  p+2,1\right)  \frac{x}{\left[  n\right]
^{p+1}}+\left(  \mathbb{S}_{q}\left(  p+2,2\right)  -2\mathbb{S}_{q}\left(
p+2,1\right)  \right)  \frac{x^{2}}{\left[  n\right]  ^{p}}\right) \\
&  =w_{p}\left(  x\right)  \frac{x}{\left[  n\right]  }\mathcal{P}_{p}\left(
q;x\right)  ,
\end{align*}
where $\mathcal{P}_{p}\left(  q;x\right)  $ is a polynomial of degree $p$.
Therefore one has%
\[
w_{p}\left(  x\right)  M_{n,q}\left(  \left(  t-x\right)  ^{2}t^{p};x\right)
\leq K_{2}\left(  q,p\right)  \frac{x}{\left[  n\right]  }.
\]

\end{proof}

Our first main result in this section is a local approximation property of
$M_{n,q}$ stated below.

\begin{theorem}
\label{t:loc}There exists an absolute constant $C>0$ such that%
\[
w_{p}\left(  x\right)  \left\vert M_{n,q}\left(  g;x\right)  -g\left(
x\right)  \right\vert \leq K_{3}\left(  q,p\right)  \left\Vert g^{\prime
\prime}\right\Vert \frac{x}{\left[  n\right]  },
\]
where $g\in C_{p}^{2}$, $q>1$ and $x\in\lbrack0,\infty)$.
\end{theorem}

\begin{proof}
Using the Taylor formula%
\[
g\left(  t\right)  =g\left(  x\right)  +g^{\prime}\left(  x\right)  )\left(
t-x\right)  +\int_{x}^{t}\int_{x}^{s}g^{\prime\prime}\left(  u\right)
du\ ds,\ \ \ \ \ \ g\in C_{p}^{2},
\]
we obtain that%
\begin{align*}
w_{p}\left(  x\right)  \left\vert M_{n,q}\left(  g;x\right)  -g\left(
x\right)  \right\vert  &  =w_{p}\left(  x\right)  \left\vert M_{n,q}\left(
\int_{x}^{t}\int_{x}^{s}g^{\prime\prime}\left(  u\right)  du\ ds;x\right)
\right\vert \\
&  \leq w_{p}\left(  x\right)  M_{n,q}\left(  \left\vert \int_{x}^{t}\int
_{x}^{s}g^{\prime\prime}\left(  u\right)  du\ ds\right\vert ;x\right) \\
&  \leq w_{p}\left(  x\right)  M_{n,q}\left(  \left\Vert g^{\prime\prime
}\right\Vert _{p}\left\vert \int_{x}^{t}\int_{x}^{s}\left(  1+u^{m}\right)
du\ ds\right\vert ;x\right) \\
&  \leq w_{p}\left(  x\right)  \frac{1}{2}\left\Vert g^{\prime\prime
}\right\Vert _{p}M_{n,q}\left(  \left(  t-x\right)  ^{2}\left(  1/w_{p}\left(
x\right)  +1/w_{p}\left(  t\right)  \right)  ;x\right) \\
&  \leq\frac{1}{2}\left\Vert g^{\prime\prime}\right\Vert _{p}\left(
M_{n,q}\left(  \left(  t-x\right)  ^{2};x\right)  +w_{p}\left(  x\right)
M_{n,q}\left(  \left(  t-x\right)  ^{2}w_{p}\left(  t\right)  ;x\right)
\right) \\
&  \leq K_{3}\left(  q,x\right)  \left\Vert g^{\prime\prime}\right\Vert
_{p}\frac{x}{\left[  n\right]  }.
\end{align*}

\end{proof}

Now we consider the modified Steklov means%
\[
f_{h}(x):=\frac{4}{h^{2}}\int\limits_{0}^{\frac{h}{2}}\int\limits_{0}%
^{\frac{h}{2}}\left[  2f(x+s+t)-f(x+2(s+t))\right]  dsdt.
\]
$f_{h}(x)$ has the following properties:%
\[
f(x)-f_{h}(x)=\frac{4}{h^{2}}\int\limits_{0}^{\frac{h}{2}}\int\limits_{0}%
^{\frac{h}{2}}\Delta_{s+t}^{2}f(x)dsdt,\ \ f_{h}^{\prime\prime}(x)=h^{-2}%
\left(  8\Delta_{\frac{h}{2}}^{2}f(x)-\Delta_{h}^{2}f(x)\right)
\]
and therefore%
\[
\left\Vert f-f_{h}\right\Vert _{p}\leq\omega_{p}^{2}(f;h),\ \ \ \left\Vert
f_{h}^{\prime\prime}\right\Vert _{p}\leq\frac{1}{9h^{2}}\omega_{p}^{2}(f;h).
\]
We have the following direct approximation theorem:

\begin{theorem}
\label{t:global}For every $p\in\mathbb{N\cup}\left\{  0\right\}  ,f\in C_{p}$
and $x\in\lbrack0,\infty),\ q>1,$ we have%
\[
w_{p}(x)\left\vert M_{n,q}\left(  f;x\right)  -f(x)\right\vert \leq
M_{p}\omega_{p}^{2}\left(  f;\sqrt{\dfrac{x}{[n]}}\right)  =M_{p}\omega
_{p}^{2}\left(  f;\sqrt{\dfrac{\left(  q-1\right)  x}{\left(  q^{n}-1\right)
}}\right)  .
\]
Particularly, if $Lip_{p}^{2}\alpha$ for some $\alpha\in(0,2],$ then%
\[
w_{p}(x)\left\vert M_{n,q}\left(  f;x\right)  -f(x)\right\vert \leq
M_{p}\left(  \dfrac{x}{[n]}\right)  ^{\frac{\alpha}{2}}%
\]

\end{theorem}

\begin{proof}
For $f\in C_{p}$ and $h>0$%
\[
\left\vert M_{n,q}(f;x)-f(x)\right\vert \leq\left\vert M_{n,q}\left(  \left(
f-f_{h}\right)  ;x\right)  -(f-f_{h})(x)\right\vert +\left\vert M_{n,q}\left(
f_{h};x\right)  -f_{h}(x)\right\vert
\]
and therefore
\begin{align*}
w_{p}(x)\left\vert M_{n,q}\left(  f;x\right)  -f(x)\right\vert  &
\leq\left\Vert f-f_{h}\right\Vert _{p}\left(  w_{p}(x)M_{n,q}\left(  \frac
{1}{w_{p}(t)};x\right)  +1\right) \\
&  +K_{3}\left(  q,p\right)  \left\Vert f_{h}^{\prime\prime}\right\Vert
_{p}\frac{x}{\left[  n\right]  }.
\end{align*}
Since $w_{p}(x)M_{n,q}\left(  \frac{1}{w_{p}(t)};x\right)  \leq K_{1}\left(
q,p\right)  $, we get that%
\[
w_{p}(x)\left\vert M_{n,q}\left(  f;x\right)  -f(x)\right\vert \leq M\left(
q,p\right)  \omega_{p}^{2}(f;h)\left[  1+\frac{x}{\left[  n\right]  h^{2}%
}\right]
\]
Thus, choosing $h=\sqrt{\dfrac{x}{[n]}}$, the proof is completed.
\end{proof}

\begin{corollary}
\label{c:convergence}If $p\in\mathbb{N\cup}\left\{  0\right\}  ,f\in C_{p},$
$q>1$ and $x\in\lbrack0,\infty),\ $then%
\[
\lim_{n\rightarrow\infty}M_{n,q}\left(  f;x\right)  =f\left(  x\right)  .
\]
This converegnce is uniform on every $\left[  a,b\right]  $, $0\leq a<b.$
\end{corollary}

\begin{remark}
Theorem \ref{t:global} shows the rate of approximation by the $q$-Szasz
operators ($q>1$) is of order $q^{-n}$ versus $1/n$ for the classical
Szasz-Mirakjan operators.
\end{remark}

\section{Convergence of $q$-Sz\'{a}sz operators}

An interesting problem is to determine the class of all continuous functions
$f$ such that $M_{n,q}\left(  f\right)  $ converges to $f$ uniformly on the
whole interval $\left[  0,\infty\right)  $ as $n\rightarrow\infty.$ This
problem was investigated by Totik \cite[Theorem 1]{totik} and de la Cal
\cite[Theorem1]{cal}. The following result is a $q$-analogue of Theorem 1
\cite{cal}.

\begin{theorem}
Assume that $f:\left[  0,\infty\right)  \rightarrow R$ is bounded or uniformly
continuus. Let
\[
f^{\ast}\left(  z\right)  =f\left(  z^{2}\right)  ,\ \ z\in\left[
0,\infty\right)  .
\]
We have, for all $t>0$ and $x\geq0$,%
\begin{equation}
\left\vert M_{n,q}\left(  f;x\right)  -f\left(  x\right)  \right\vert
\leq2\omega\left(  f^{\ast};\sqrt{\frac{1}{\left[  n\right]  }}\right)  .
\label{s1}%
\end{equation}
Therefore, $M_{n,q}\left(  f;x\right)  $ converges to $f$ uniformly on
$\left[  0,\infty\right)  $ as $n\rightarrow\infty$, whenever $f^{\ast}$ is
uniformly continuous.
\end{theorem}

\begin{proof}
By the definition of $f^{\ast}$ we have
\[
M_{n,q}\left(  f;x\right)  =M_{n,q}\left(  f^{\ast}\left(  \sqrt{\cdot
}\right)  ;x\right)  .
\]
Thus we can write%
\begin{align*}
\left\vert M_{n,q}\left(  f;x\right)  -f\left(  x\right)  \right\vert  &
=\left\vert M_{n,q}\left(  f^{\ast}\left(  \sqrt{\cdot}\right)  ;x\right)
-f^{\ast}\left(  \sqrt{x}\right)  \right\vert \\
&  =\left\vert \sum_{k=0}^{\infty}\left(  f^{\ast}\left(  \sqrt{\frac{\left[
k\right]  }{\left[  n\right]  }}\right)  -f^{\ast}\left(  \sqrt{x}\right)
\right)  s_{n,k}\left(  q;x\right)  \right\vert \\
&  \leq\sum_{k=0}^{\infty}\left\vert \left(  f^{\ast}\left(  \sqrt
{\frac{\left[  k\right]  }{\left[  n\right]  }}\right)  -f^{\ast}\left(
\sqrt{x}\right)  \right)  \right\vert s_{n,k}\left(  q;x\right) \\
&  \leq\sum_{k=0}^{\infty}\omega\left(  f^{\ast};\left\vert \sqrt
{\frac{\left[  k\right]  }{\left[  n\right]  }}-\sqrt{x}\right\vert \right)
s_{n,k}\left(  q;x\right) \\
&  \leq\sum_{k=0}^{\infty}\omega\left(  f^{\ast};\frac{\left\vert \sqrt
{\frac{\left[  k\right]  }{\left[  n\right]  }}-\sqrt{x}\right\vert }%
{M_{n,q}\left(  \left\vert \sqrt{\cdot}-\sqrt{x}\right\vert ;x\right)
}M_{n,q}\left(  \left\vert \sqrt{\cdot}-\sqrt{x}\right\vert ;x\right)
\right)  s_{n,k}\left(  q;x\right)  .
\end{align*}
Finally, from the inequality%
\[
\omega\left(  f^{\ast};\alpha\delta\right)  \leq\left(  1+\alpha\right)
\omega\left(  f^{\ast};\delta\right)  ,\ \ \ \ \ \alpha,\delta\geq0,
\]
we obtain%
\begin{align*}
\left\vert M_{n,q}\left(  f;x\right)  -f\left(  x\right)  \right\vert  &
\leq\omega\left(  f^{\ast};M_{n,q}\left(  \left\vert \sqrt{\cdot}-\sqrt
{x}\right\vert ;x\right)  \right)  \sum_{k=0}^{\infty}\left(  1+\frac
{\left\vert \sqrt{\frac{\left[  k\right]  }{\left[  n\right]  }}-\sqrt
{x}\right\vert }{M_{n,q}\left(  \left\vert \sqrt{\cdot}-\sqrt{x}\right\vert
;x\right)  }\right)  s_{n,k}\left(  q;x\right) \\
&  =2\omega\left(  f^{\ast};M_{n,q}\left(  \left\vert \sqrt{\cdot}-\sqrt
{x}\right\vert ;x\right)  \right)  .
\end{align*}
In order to complete the proof we need to show that we have for all $t>0$ and
$x>0$,%
\[
M_{n,q}\left(  \left\vert \sqrt{\cdot}-\sqrt{x}\right\vert ;x\right)
\leq\sqrt{\frac{1}{\left[  n\right]  }}.
\]
Indeed we obtain from the Cauchy-Schwarz inequality%
\begin{align*}
M_{n,q}\left(  \left\vert \sqrt{\cdot}-\sqrt{x}\right\vert ;x\right)   &
=\sum_{k=0}^{\infty}\left\vert \sqrt{\frac{\left[  k\right]  }{\left[
n\right]  }}-\sqrt{x}\right\vert s_{n,k}\left(  q;x\right) \\
&  =\sum_{k=0}^{\infty}\frac{\left\vert \frac{\left[  k\right]  }{\left[
n\right]  }-x\right\vert }{\sqrt{\frac{\left[  k\right]  }{\left[  n\right]
}}+\sqrt{x}}s_{n,k}\left(  q;x\right)  \leq\frac{1}{\sqrt{x}}\sum
_{k=0}^{\infty}\left\vert \frac{\left[  k\right]  }{\left[  n\right]
}-x\right\vert s_{n,k}\left(  q;x\right) \\
&  \leq\frac{1}{\sqrt{x}}\sqrt{\sum_{k=0}^{\infty}\left\vert \frac{\left[
k\right]  }{\left[  n\right]  }-x\right\vert ^{2}s_{n,k}\left(  q;x\right)
}=\frac{1}{\sqrt{x}}\sqrt{M_{n,q}\left(  \left(  \cdot-x\right)
^{2};x\right)  }\\
&  =\frac{1}{\sqrt{x}}\sqrt{\frac{1}{\left[  n\right]  }x}=\sqrt{\frac
{1}{\left[  n\right]  }}%
\end{align*}
showing (\ref{s1}), and completing the proof.
\end{proof}

Next we prove Voronovskaja type result for $q$-Sz\'{a}sz-Mirakjan operators.

\begin{theorem}
Assume that $q\in\left(  1,\infty\right)  $. For any $f\in C_{p}^{2}$ the
following equality holds%
\[
\lim_{n\rightarrow\infty}\left[  n\right]  \left(  M_{n,q}\left(  f;x\right)
-f\left(  x\right)  \right)  =\frac{1}{2}f^{\prime\prime}\left(  x\right)  x
\]
for every $x\in\left[  0,\infty\right)  $.
\end{theorem}

\begin{proof}
Let $x\in\left[  0,\infty\right)  $ be fixed. By the Taylor formula we may
write%
\begin{equation}
f\left(  t\right)  =f\left(  x\right)  +f^{\prime}\left(  x\right)  \left(
t-x\right)  +\frac{1}{2}f^{\prime\prime}\left(  x\right)  \left(  t-x\right)
^{2}+r\left(  t;x\right)  \left(  t-x\right)  ^{2}, \label{v3}%
\end{equation}
where $r\left(  t;x\right)  $ is the Peano form of the remainder, $r\left(
\cdot;x\right)  \in C_{p}$ and $\lim_{t\rightarrow x}r\left(  t;x\right)  =0$.
Applying $M_{n,q}$ to (\ref{v3}) we obtain
\begin{multline*}
\left[  n\right]  \left(  M_{n,q}\left(  f;x\right)  -f\left(  x\right)
\right)  =f^{\prime}\left(  x\right)  \left[  n\right]  M_{n,q}\left(
t-x;x\right) \\
+\frac{1}{2}f^{\prime\prime}\left(  x\right)  \left[  n\right]  M_{n,q}\left(
\left(  t-x\right)  ^{2};x\right)  +\left[  n\right]  M_{n,q}\left(  r\left(
t;x\right)  \left(  t-x\right)  ^{2};x\right)  .
\end{multline*}
By the Cauchy-Schwartz inequality, we have%
\begin{equation}
M_{n,q}\left(  r\left(  t;x\right)  \left(  t-x\right)  ^{2};x\right)
\leq\sqrt{M_{n,q}\left(  r^{2}\left(  t;x\right)  ;x\right)  }\sqrt
{M_{n,q}\left(  \left(  t-x\right)  ^{4};x\right)  }. \label{v4}%
\end{equation}
Observe that $r^{2}\left(  x;x\right)  =0.$ Then it follows from Corollary
\ref{c:convergence} that%
\begin{equation}
\lim_{n\rightarrow\infty}M_{n,q}\left(  r^{2}\left(  t;x\right)  ;x\right)
=r^{2}\left(  x;x\right)  =0. \label{v5}%
\end{equation}
Now from (\ref{v4}), (\ref{v5}) and Lemma \ref{lema} we get immediately%
\[
\lim_{n\rightarrow\infty}\left[  n\right]  M_{n,q}\left(  r\left(  t;x\right)
\left(  t-x\right)  ^{2};x\right)  =0.
\]
The proof is completed.
\end{proof}


\bigskip

\begin{thebibliography}{99}                                                                                               %


\bibitem {altomare}F. Altomare, M. Campiti, Korovkin-type approximation theory
and its applications. de Gruyter Studies in Mathematics, 17. Walter de Gruyter
\& Co., Berlin, 1994. xii+627 pp.

\bibitem {lupas}A. Lupa\c{s}, A $q$-analogue of the Bernstein operator.
Seminar on numerical and statistical calculus, University of Cluj-Napoca 9:
85-92 (1987 )

\bibitem {phil}G.M. Phillips, Bernstein polynomials based on the $q$-integers,
Ann. Numer. Math., 4(1997) 511-518.

\bibitem {sof1}A. II'inskii, S.Ostrovska, Convergence of generalized Bernstein
polynomials, J. Approx. Theory, 116(1)(2002) 100-112.

\bibitem {sof2}S. Ostrovska, $q$-Bernstein polynomials and their iterates, J.
Approx. Theory, 123(2)(2003) 232-255.

\bibitem {sof3}S. Ostrovska, On the limit $q$-Bernstein operators, Math.
Balkanica (N.S), 18(2004) 165-172.

\bibitem {sof4}S. Ostrovska, On the Lupa\c{s} $q$-analogue of the Bernstein
operator, \textit{Rocky Mountain Journal of Mathematics}, \textbf{36},
1615-1629 (2006).

\bibitem {mah}N.I. Mahmudov, Korovkin-type Theorems and Applications, Central
European Journal of Mathematics, DOI: 10.2478/s11533-009-0006-7.

\bibitem {vid}V.S. Videnskii, On some classes of $q$-parametric positive
operators, Operator Theory: Advances and Applications, Birkhauser, Basel
158(2005) 213-222..

\bibitem {hep2}H. Wang, F. Meng, The rate of convergence of $q$-Bernstein
polynomials for $0<q<1$, J. Approx. Theory 136(2)(2005) 151-158.

\bibitem {trif}T. Trif, M eyer-Konig and Zeller operators based on the
q-integers. Rev . Anal. Numer. Theor. Approx. 29: 221-229 (2000).

\bibitem {dogru}O . Dogru, O. Duman, Statistical approximation of Meyer-Konig
and Zeller operators based on q-integers. Publ. Math. Debrecen 68 :199 -214 (2006).

\bibitem {hep3}H. Wang, Properties of convergence for the $q$-Meyer-K\"{o}nig
and Zeller operators. J. Math. Anal. Appl. 335 (2007), no. 2, 1360--1373..

\bibitem {aral2}A. Aral, O. Do\u{g}ru, Bleimann, Butzer, and Hahn operators
based on the $q$-integers. \textit{J. Inequal. Appl.} 2007, Art. ID 79410.

\bibitem {mahsab}N. I. Mahmudov and P. Sabanc\i gil, $q$-Parametric Bleimann
Butzer and Hahn Operators, \textit{Journal of Inequalities and Applications,
}\textbf{2008} (2008), Article ID 816367, 15 pages do:10.1155/2008/816367.

\bibitem {mahsab2}Mahmudov, N. I.; Sabancigil, P., A $q$-analogue of the
Meyer-K\"{o}nig and Zeller operators, Bull. Malaysian Math. Soc. Sci.., In press

\bibitem {gupta}A. Aral, V. Gupta, The $q$-derivative and applications to
$q$-Sz\'{a}sz Mirakyan operators. Calcolo 43 (2006), no. 3, 151--170.

\bibitem {aral}A. Aral, A generalization of Sz\'{a}sz--Mirakyan operators
based on $q$-integers, Mathematical and Computer Modelling 47 (2008), no.
9-10, 1052-1062.

\bibitem {mah2}N.I. Mahmudov, On $q$-parametric Sz\'{a}sz-Mirakjan operators,
Med. J. Math., In press.

\bibitem {radu}C. Radu, On statistical approximation of a general class of
positive linear operators extended in q-calculus, Appl. Math. Comput. (2009),
doi: 10.1016/j.amc.2009.08.023.

\bibitem {ag}A. Aral and V. Gupta, On certain $q$-Baskakov operators, prepeint.

\bibitem {becker}M. Becker: Global approximation theorems for
Sz\'{a}sz-Mirakjan and Baskakov operators in polynomial weight spaces, Indiana
Univ. Math. J., 27(1)(1978), 127-142.

\bibitem {cal}J. de la Cal and J. C\'{a}rcamo, On uniform approximation by
some classical Bernstein-type operators, J. Math. Anal. Appl. 279 (2003) 625--638.

\bibitem {ernst}T. Ernst, The history of $q$-calculus and a new method .
U.U.D.M . Report 2000, 16 , Uppsala, Departament of Mathematics, Uppsala
University 2000.

\bibitem {kac}V. Kac and P. Cheung, Quantum Calculus, Universitext,
Springer-Verlag, New York, 2002.

\bibitem {szasz}O. Sz\'{a}sz, Generalization of S. Bernstein's polynomials to
the infinite interval. J. Research Nat. Bur. Standards 45, 239--245 (1950).

\bibitem {totik}V. Totik, Uniform approximation by Sz\'{a}sz--Mirakjan type
operators, Acta Math. Hungar. 41 (1983) 291--307.
\end{thebibliography}
\end{document}